\documentclass[12pt,twoside]{article}
\title{On
$n$-Perfect Rings  and  Cotorsion Dimension}
\date{}

\usepackage{amsfonts}
\usepackage{amsmath}
\usepackage{amssymb}
\usepackage{latexsym}
\usepackage{graphicx}

\newtheorem{thm}{\bf Theorem}[section]

\newtheorem{lem}[thm]{\bf Lemma}
\newtheorem{prop}[thm]{\bf Proposition}
\newtheorem{defn}[thm]{\bf Definition}

\newtheorem{rem}[thm]{\bf Remark}

\newtheorem{exmp}[thm]{\bf Example}

\catcode`\ç=13
\defç{\c{c}}
\catcode`\é=13
\defé{\'e}
\catcode`\à=13
\defà{\`a}
\catcode`\è=13
\defè{\`e}
\catcode`\â=13
\defâ{\^a}
\catcode`\ù=13
\defù{\`u}
\catcode`\ê=13
\defê{\^e}
\catcode`\î=13
\defî{\^\i}
\catcode`\ô=13
\defô{\^o}
\newcommand{\field}[1]{\mathbb{#1}}
\newcommand{\C}{\field{C}}

\newcommand{\Q }{\field{Q}}
\newcommand{\Z }{\field{Z}}
\newcommand{\N }{\field{N}}

\def\proof{{\parindent0pt {\bf Proof.\ }}}

\def\wdim{{\rm wdim}}
\def\gldim{{\rm gldim}}

\def\cot.D{{\rm C\!-\!gldim}}

\def\FPD{{\rm FPD}}
\def\FID{{\rm FID}}
\def\FFD{{\rm FFD}}
\def\FCD{{\rm FCD}}

\def\pd{{\rm pd}}
\def\fd{{\rm fd}}
\def\id{{\rm id}}

\def\cd{{\rm cd}}

\def\Ext{{\rm Ext}}
\def\Tor{{\rm Tor}}
\def\Hom{{\rm Hom}}

\def\sup{{\rm sup}}

\def\qf{{\rm qf}}

\newcommand{\cqfd}
{\hspace{1cm}
\rule{2mm}{2mm}%
\medbreak%
\par%
}

\begin{document}
\thispagestyle{empty}
\maketitle \vspace*{-1.5cm}
\begin{center}{\large\bf Driss Bennis and Najib Mahdou}

\bigskip

\small{Department of Mathematics, Faculty of Science and
Technology of Fez,\\ Box 2202, University S. M.
Ben Abdellah Fez, Morocco\\[0.12cm]
driss\_bennis@hotmail.com\\
mahdou@hotmail.com}
\end{center}

\bigskip\bigskip
\noindent{\large\bf Abstract.} A ring is called $n$-perfect
($n\geq 0$), if every flat module has projective dimension less or
equal than $n$.\\
In this paper, we show that the $n$-perfectness relates, via
homological approach, some homological dimensions of rings. We study
$n$-perfectness in some known ring's construction. Finally, several
examples of $n$-perfect rings satisfying special conditions are
given.\bigskip

\medskip



\begin{section}{Introduction} Throughout this paper  all   rings
are commutative with identity element and all modules are unitary.
For a ring $R$ and an $R$-module $M$, we use $\pd_R(M),\ \id_R(M)$,
and $\fd_R(M)$ to denote, respectively, the classical projective,
injective, and flat dimension of $M$. We use $\gldim(R)$ and
$\wdim(R)$ to denote, respectively, the classical global and weak
dimension of $R$. If $R$ is integral, we denote its quotient field
by $\qf(R)$.\bigskip

In \cite{bass}, Bass proved that the perfect rings are those rings
such that every flat module is projective. He links these rings with
the finitistic projective dimension of rings. Recall that the
finitistic projective dimension of a ring  $R$, denoted  by
$\FPD(R)$, is defined by: $\FPD(R)=\sup\{\pd_R(M)|M\ R\!-\!module\
with\ \pd_R(M)<\infty\}.$\\
From  \cite[Kaplansky's Theorem, page 466]{bass} and \cite[Part I,
3. Example (6)]{bass}, the following are equivalent, for a
commutative ring  $R$:
\begin{enumerate}
    \item $R$ is perfect;
    \item $R$ is a finite direct product of local (not necessarily Noetherian) rings, each with
$T$-nilpotent  maximal ideal (i.e., if we pick a sequence $a_1,
a_2,...$ of elements in the maximal ideal, then for some index
$m$, $a_1a_2...a_m=0$);
    \item   $\FPD(R)=0$.
\end{enumerate}

\indent Later, in \cite[Proposition 6]{Jensen}, Jensen proved that,
for a ring $R$, if $\FPD(R)=n$ ($n\geq 0$), then every flat
$R$-module has projective dimension at most $n$. In several
situations the Jensen's result implies the desired properties. But,
unfortunately, the finitistic projective dimension is not known to
be finite only in few cases (for example for Noetherian rings with
finite Krull dimension \cite[Theorem 3.2.6]{RG}), and its finiteness
remains an open problem. Then, it seems appropriate to investigate
the rings over which every flat module has projective dimension at
most $n$, where $n$ is a fixed positive integer. In \cite[Definition
1.1]{n-perfect}, Enochs, Jenda, and L\'{o}pez-Ramos called these
rings $n$-perfect. These rings were homologically characterized by
the cotorsion dimension introduced by Ding and Mao as follows:

\begin{defn}[\cite{DM}]\label{def-Cot-dim} \textnormal{Let $R$ be a ring.\\
\indent The cotorsion dimension of an $R$-module $M$, denoted by
$\cd_R(M)$, is the least positive integer $n$ for which
$\Ext^{n+1}_{R}(F, C)=0$ for all flat $R$-modules
$F$.\\
\indent The  global cotorsion dimension of $R$, denoted by
$\cot.D(R)$, is the quantity: $\cot.D(R)=\sup\{\cd_R(M)|\ M \
R\!-\!module\}$.}
\end{defn}
Namely, the modules of cotorsion dimension $0$ are the known
cotorsion modules (see \cite[Definition 3.1.1]{Xu}). We have:

\begin{prop}[\cite{DM}, Theorem 19.2.5(1)]\label{pro-cot=nperfect}
 For a positive integer $n$, $R$ is
$n$-perfect  if and only if $\cot.D(R)\leq n$.
\end{prop}

So, Jensen's result above can be written again as follows:

\begin{prop}[\cite{Jensen}, Proposition 6]\label{pro-Jensen}
For any ring $R$, $\cot.D(R)\leq \FPD(R).$
\end{prop}

Early, in \cite[page 87]{RG}, Gruson and Raynaud defined $d(R)$ as
the supremum of the projective dimensions of  all flat $R$-modules.
Then, $d(R)$ coincides with the  global cotorsion dimension of $R$.
They shortly studied this invariant of rings and they mentioned that
Jensen had an example of a ring $R$ that satisfies the strict
inequality $\cot.D(R)< \FPD(R).$ One of our aims in this paper is to
give further concrete examples of rings with the strict inequality,
which means that  the converse of Jensen's result is not true in
general (see Section 4).\bigskip

In Section 2, we give some general results on the  global cotorsion
dimension of rings. Mainly, we extent the inequality $ \gldim(R)\leq
\cot.D(R)+ \wdim(R)$ established by Ding and Mao in \cite[Theorem
19.2.14]{DM} to the finitistic dimensions (see
Theorems \ref{thm-C-FPD-FFD} and  \ref{thm-wdim-FID-FCD}).\\
In Section 3, we investigate $n$-perfectness in some known ring's
constructions, such that we compute the  global cotorsion dimension
of polynomial rings, finite direct products of rings, and $D+M$
rings. This study allows us to give various examples of $n$-perfect
rings satisfying special conditions. This is given  in Section 4.

\end{section}
\bigskip\bigskip


\begin{section}{General results}
In \cite[Theorem 19.2.14]{DM}, the  global cotorsion dimension of
rings is used to relate the weak and the global dimensions as
follows:

\begin{thm}[\cite{DM}, Theorem 19.2.14]\label{thm-C-gldim-wdim}
For any ring $R$, the following inequality $
  \gldim(R)\leq \cot.D(R)+ \wdim(R) $ holds true.\\
In particular: \begin{itemize}
    \item If $\cot.D(R)=0$ (i.e., $R$ is perfect), then $\wdim(R)=
\gldim(R)$.
    \item If $\wdim(R)=0$ (i.e., $R$ is von Neumann regular), then $\cot.D(R)=
\gldim(R)$.
\end{itemize}
\end{thm}

The next theorem generalizes this result to the finitistic
projective and flat dimensions. Recall that the finitistic flat
dimension of $R$, $\FFD(R)$, is defined as follows:
$\FFD(R)=\sup\{\fd_R(M)|M\ R\!-\!module\ with\ \fd_R(M)<\infty\}$.

\begin{thm}\label{thm-C-FPD-FFD}
For any ring $R$, the following inequalities  $\FFD(R)\leq
\FPD(R)\leq\cot.D(R) +  \FFD(R) $  hold true.
\end{thm}
\proof To prove the inequality  $\FFD(R)\leq \FPD(R)$ we can assume
that  $\FPD(R)=n$ is finite. Consider an $R$-module $M$ with finite
flat dimension. From Jensen's Proposition \ref{pro-Jensen}, $M$ has
also finite projective dimension, which is at most $n$. Then, we
have $\fd_R(M)\leq \pd_R(M)\leq n$. This means
that $\FFD(R)\leq \FPD(R)$, as desired.\\
Now we prove the second inequality. For that we can assume that
$\cot.D(R)=n$ and $\FFD(R)=m$ are finite. Consider an $R$-module $M$
with finite projective dimension, then it has finite flat dimension,
which is at most $m$. Then, there exists an exact sequence of
$R$-modules $$0\rightarrow F \rightarrow P_{m-1}\rightarrow \cdots
\rightarrow P_{0} \rightarrow  N \rightarrow 0,$$ where $P_i$ are
projective  and   $F$ is flat. From Proposition
\ref{pro-cot=nperfect}, $\pd_R(F)\leq n$. Finally, using the above
sequence, we get $\pd_R(M)\leq n+m$. This completes the
proof.\cqfd\bigskip

In Section 4, Example \ref{exmp3} shows that the upper bound on the
finitistic projective dimension in Theorem \ref{thm-C-FPD-FFD}
above is the best possible upper bound.\\
As mentioned in the introduction, in \cite[page  87]{RG}, $d(R)$
denotes, for a ring $R$, the supremum of the projective dimensions
of all flat $R$-modules, then it coincides with the  global
cotorsion dimension of $R$. The authors mentioned that Jensen
constructed an example showing that the inequality $\cot.D(R)\leq
\FPD(R)$ can be strict. In Section 4, we give several examples
showing that each of the inequalities of both Theorems
\ref{thm-C-gldim-wdim} and \ref{thm-C-FPD-FFD}
 can be strict.\bigskip

The following gives a situation in which the upper bound on the
finitistic projective dimension (and then of the global dimension)
in Theorem  \ref{thm-C-FPD-FFD} (and Theorem \ref{thm-C-gldim-wdim})
is reached. Later, Example  \ref{exmp3} (and Example \ref{exmp1})
provides a ring  satisfying this situation.

\begin{prop}\label{pro-FPD-eqUp}
If a ring $R$ satisfies $\FFD(R)=1$ and
$\cot.D(R)<\FPD(R)<\infty$, then: $ \FPD(R)=\cot.D(R) +1.$
\end{prop}
\proof Let $\FPD(R)=n<\infty$ for an integer $n\geq 1$. Then, there
is an $R$-module  $M$ that satisfies $\pd_R(M)=n$. Hence, for a
short exact sequence of $R$-modules $0\rightarrow F \rightarrow P
\rightarrow M \rightarrow 0$, where $P$ is   projective  and then
$F$ is  flat (since $\FFD(R)=1$),  we have $\pd_R(F)=n-1$. Then, the
desired equality holds since $\cot.D(R)<\FPD(R)$.\cqfd\bigskip

Recall that the finitistic injective dimension of a ring $R$,
denoted by $\FID(R)$, is defined by: $\FID(R)=\sup\{\id_R(M)|M\
R\!-\!module\ with\ \id_R(M)<\infty\}.$ Similarly, we can define the
finitistic cotorsion dimension of a ring  $R$, denoted by $\FCD(R)$,
as follows:
$\FCD(R)=\sup\{\cd_R(M)|M\ R\!-\!module\ with\ \cd_R(M)<\infty\}.$\\
The following result is another extension of Theorem
\ref{thm-C-gldim-wdim}.

\begin{thm}\label{thm-wdim-FID-FCD}
For any ring $R$ with finite weak dimension, the following
inequalities $ \FCD(R)\leq\FID(R) \leq \FCD(R) + \wdim(R) $ hold
true.
\end{thm}
The proof  involves the following lemma which relates the cotorsion
dimension and the injective dimension of modules.

\begin{lem}\label{lem-cd-id}
Let $R$ be a ring. For any  $R$-module  $M$, the following
inequalities
 $\cd(M)\leq \id(M)\leq \cd(M)+\wdim(R) $  hold
true.
\end{lem}
\proof The proof of the first inequality is easy.\\
To prove the second inequality, we can assume that $ \cd(M)=m$ and
$\wdim(R)=n$ are finite. Let $N$ be any $R$-module, and consider an
exact sequence $0\rightarrow F \rightarrow P_{n-1}\rightarrow \cdots
\rightarrow P_{0} \rightarrow  N \rightarrow 0$, where $P_i$ are
projective modules and then $F$ is a flat module (since
$\wdim(R)=n$). We have $\Ext^{m+n+1}(N,M)\cong\Ext^{m+1}(F,M)=0$
(since $ \cd(M)=m$). Therefore, $\id(M)\leq m+n$, as
desired.\cqfd\bigskip

\noindent\textbf{Proof of Theorem \ref{thm-wdim-FID-FCD}.} We prove
the inequality $ \FCD(R)\leq\FID(R)$.  We can suppose that
$\FID(R)=n$ is finite. Let $M$ be an $R$-module with finite
cotorsion dimension. From Lemma \ref{lem-cd-id}, $\id(M)\leq
\cd(M)+\wdim(R)$  which is finite.
Then, $\id(M)$ is finite and so  $\cd(M)\leq \id(M)\leq n$, as desired.\\
We prove the inequality $ \FID(R) \leq \FCD(R)+ \wdim(R).$ We can
suppose that $\wdim(R) =n$ and  $\FCD(R)=m$ are finite. Let $M$ be
an $R$-module with finite injective dimension. Then, from Lemma
\ref{lem-cd-id}, $\id(M)\leq \cd(M)+\wdim(R)\leq m+n$. This implies
the desired inequality.\cqfd\bigskip

In \cite[Remark 3.3.3(1), page  87]{RG}, the authors ask whether,
for an integral Noetherian domain $R$ with quotient field $Q$,
$\pd_R(Q)=\cot.D(R)$. Next, we analyze  some situations in which we
have an affirmative answer. But, for non-Noetherian domains, Example
\ref{exmp2} shows that we have a negative answer (see Remark
\ref{rem-dom}).

\begin{prop}\label{pro-C-gldim-domain}
Let $R$ be a  domain which is not a field, and let $Q$ be the
quotient field of $R$. Then, the following assertions hold true:
\begin{enumerate}
    \item If $\cot.D(R)\leq 1$, then $\cot.D(R)= \pd_R(Q)$.
    \item If $\cot.D(R)= \FPD(R)<\infty$, then $\cot.D(R)= \pd_R(Q)$.
    \item If $\gldim(R)=\pd_R(Q)+1$, then  $\cot.D(R)=\pd_R(Q)$.
\end{enumerate}
\end{prop}
\proof 1. Since the quotient field  of a domain which is not a field
is  flat but  not  projective, we  deduce that:
\begin{itemize}
    \item If $\cot.D(R)=0$ (i.e., $R$ is perfect), then $R$ is a
    field. So, $R=Q$  and $\cot.D(R)= \pd_R(Q)=0$.
    \item  If $\cot.D(R)=1$, then $R\subsetneq Q$ and then
    $\pd_R(Q)=1$, as desired.
\end{itemize}
2. From (1), we can assume that   $0<\cot.D(R)= \FPD(R)=n<\infty$
($n\geq 1$). Let $F$ be a flat $R$-module such that $\pd_R(F)=n$.
From \cite[Lemma 4.31 and the Remark bellow Theorem 4.33]{Rot},
there exist an index set $I$ and an $R$-module $M$ such that the
sequence $0 \rightarrow F\rightarrow Q^{(I)} \rightarrow
M\rightarrow 0$  is exact. Then, if we suppose that $\pd_R(Q)<n$,
we obtain that $\pd_R(M)=n+1$, which gives a contradiction.\\
3. We can assume that  $\pd_R(Q)=n$ for a positive integer $n$. As
above, for every flat module, there are an index set $I$ and an
$R$-module $M$ such that the sequence $0 \rightarrow F\rightarrow
Q^{(I)} \rightarrow M\rightarrow 0$  is exact. From \cite[Corollary
2(c) and Corollary 1, page 135]{Bou}, we have $\pd(F)\leq \sup
\{\pd(Q), \pd(M)-1\}\leq n$ (since $\gldim(R)=\pd_R(Q)+1=n+1$). This
completes the proof.\cqfd\bigskip

In Section 4, we construct a ring satisfying  $(1)$ and $(2)$ of
Proposition \ref{pro-C-gldim-domain}   (see Example \ref{exmp1}).

\end{section}\bigskip\bigskip


\begin{section}{$n$-Perfectness in some known ring's constructions}
In this section, we establish some results concerning the transfer
of the $n$-perfectness in some known ring's constructions.
These results will be used to construct some desired examples.\\
We begin by computing the  global cotorsion dimension of polynomial
rings.

\begin{thm}\label{thm-Syzygy-C-gldim}
Let  $R[X_1,X_2,...,X_n]$ be the polynomial ring in $n$
indeterminates over a ring $R$. Then, for a positive integer $m$,
$R$ is $m$-perfect if and only if $R[X_1,X_2,...,X_n]$ is
$m+n$-perfect.\\
In other words, $ \cot.D(R[X_1,X_2,...,X_n])=\cot.D(R)+n.$
\end{thm}
\proof  By induction, it suffices to prove the case $n=1$. We
write $R[X_1]=R[X]$.\\
Through this proof, we use $M[X]$ to denote the $R[X]$-module
$M\otimes_R
R[X]$.\\
We prove that $ \cot.D(R[X])=\cot.D(R)+1.$\\
The inequality $\cot.D(R[X])\leq \cot.D(R)+1$ is the same as
\cite[Example
2.5]{n-perfect}.\\
Conversely, assume that $\cot.D(R[X])= n+1<\infty $. Let $F$ be a
flat $R$-module, then $F[X]$ is a flat $R[X]$-module. Then,
$\pd_{R}(F)=\pd_{R}(F[X])\leq \pd_{R[X]}(F[X])\leq n+1$. This means
that $\cot.D(R)\leq n+1$. Assume that $\cot.D(R)= n+1$. Then, there
exists, from \cite[Theorem 19.2.5(1)]{DM}, a flat $R$-module $F$
such that $\cd_R(F)=n+1$. Thus, there exists a flat $R$-module $E$
such that $\Ext_{R}^{n+1}(E,F)\not = 0$. From \cite[Example 7, page
9]{Bou}, the endomorphism $\mu:F\rightarrow F$, defined by
$\mu(f)=Xf$, is injective. Then, Rees's Theorem \cite[Theorem
9.37]{Rot} gives:
$$\Ext_{R[X]}^{n+2}(E,F[X])\cong\Ext_{R}^{n+1}(E,F)\not = 0.$$
From \cite[Exercise 9.20, page 258]{Rot},
$$\Ext_{R[X]}^{n+2}(E[X],F[X])\cong\Hom_R(R[X],\Ext_{R[X]}^{n+2}(E,F[X]))\cong
(\Ext_{R[X]}^{n+2}(E,F[X]))^{\N}\not = 0.$$ Then,
$\cd_{R[X]}(F[X])\geq n+2$, which contradicts
$\cot.D(R[X])=n+1$.\cqfd\bigskip

Now, we compute the  global cotorsion dimension of a finite direct
product of rings. Namely, we extent  the obvious fact that a finite
direct product of rings is perfect if  and only if  each of these
rings is perfect.

\begin{thm}\label{thm-product-C-gldim}
Let $\{R_{i}\}_{i=1,...,m}$  be a family of rings. Then,  for a
positive integer $n$, $\displaystyle\prod_{i=1}^m R_i$ is
$n$-perfect if and only if each of $R_i$ is
$n$-perfect.\\
In other words, $ \cot.D(\displaystyle\prod_{i=1}^m
R_i)=\sup\{\cot.D(R_i), 1\leq i \leq m\}.$
\end{thm}
\proof It follows by induction on $m$ and using \cite[Lemma 2.5
(2)]{Mah2001} and \cite[Lemma 3.7]{BM2}.\cqfd\bigskip

We end this section with a study of the $n$-perfectness in a
particular case of the $D+M$-constructions. These constructions have
been proven to be useful in solving many open problems and
conjectures in various contexts in ring theory (please see
\cite[Section 1, Chapter 5]{Glaz}; see also \cite{Dobs} and
\cite{Mah2001}).

\begin{thm}\label{thm-pullback-C-gldim}
Let  $T$ be a domain of the form $K+M$, where $K$ is a field and
$M(\not= 0)$ is a maximal ideal of $T$. Let $D$ be a proper
subring of $K$ such that $\qf(D)= K$. Then, for the subring
$R=D+M$
of $T$, we have $\cot.D(R)= \sup\{\cot.D(T),\cot.D(D)\}.$\\
In other words, $R$ is $n$-perfect ($n\geq0$) if and only if  $T$
and $D$ are $n$-perfect.
\end{thm}

The proof of the theorem involves the following results, which are
of independent interest.\bigskip

Next lemma shows that the $n$-perfectness descends from a ring to
any subring retract. Recall that, for  a ring homomorphism
$\psi:\, R\rightarrow S$, we say that $R$ is a subring retract of
$S$, if there exists a ring homomorphism $\phi: S\longrightarrow
R$ satisfying $\phi\psi =id_{|R}$. In this case,  $\psi$ is
injective and  the $R$-module $S$ contains $R$ as a direct summand
\cite[page 111]{Glaz}.

\begin{lem}\label{lem-retract-C-gldim}
If $R$ is a subring retract of a ring $S$, then:
$\cot.D(R)\leq\cot.D(S)$.
\end{lem}
\proof we can assume that $\cot.D(S)=n$ for a  positive integer $n$.
Let $F$ be a flat $R$-module and let $N$ be an $R$-module. We have
$\Tor^{R}_{p}(S,F)=0$ for all $p>0$. Then, from \cite[Proposition
4.1.3]{CE},
$$\Ext^{n+1}_{R}(F,N\otimes_R S) \cong \Ext^{n+1}_{S}(F\otimes_R
S,N\otimes_R S)=0.$$ Since $R$ is a direct  summand  of the
$R$-module $S$, $N$ is  a direct  summand  of the $R$-module
$N\otimes_R S$, and then $\Ext^{n+1}_{R}(F,N) $  is  a direct
summand  of the $R$-module $\Ext^{n+1}_{R}(F,N\otimes_R S) $.
Thus, $\Ext^{n+1}_{R}(F,N)=0 $, as desired.\cqfd

\begin{lem}\label{lem-pullbak2}
Let $U$ be a multiplicative set of a ring $R$. For the ring of
fractions $S=U^{-1}R$, we have: $\cot.D(S)\leq\cot.D(R)$.
\end{lem}
\proof We can assume that $\cot.D(R)=n$ for a  positive integer $n$.
Let $F$ be a flat $S$-module. Then, $F$ is a flat $R$-module, and so
$\pd_R(F)\leq n$. Thus, there is an exact sequence of $R$-modules
$$0\rightarrow P_n \rightarrow  \cdots \rightarrow
P_{0} \rightarrow  F \rightarrow 0,$$ where each $P_i$ is
projective. Then, the sequence  of $S$-modules
$$0\rightarrow P_n \otimes_R S\rightarrow \cdots \rightarrow P_{0}
\otimes_R S\rightarrow F\otimes_R S=F \rightarrow 0 $$ is an exact
sequence with each $P_i\otimes_R S$ is a projective $S$-module.
This means that $S$ is $n$-perfect and therefore $\cot.D(S)\leq
n$.\cqfd\bigskip

\noindent\textbf{Proof of Theorem \ref{thm-pullback-C-gldim}.}
First, note that $U^{-1}R=T$, where $U=D\setminus\{0\}$ (since
$\qf(D)=K$). Then, $\cot.D(T)\leq \cot.D(R)$ (by Lemma
\ref{lem-pullbak2}). And, from Lemma \ref{lem-retract-C-gldim},
$\cot.D(D)\leq \cot.D(R)$ (since $D$ is a subring retract of
$R$).\\
It remains to prove the inequality $ \cot.D(R) \leq
\sup\{\cot.D(T),\cot.D(D)\} $. We can assume that $\cot.D(T)\leq n$
and $\cot.D(D)\leq n$  for a  positive integer $n$. Let $F$ be a
flat $R$-module. Consider the following exact sequence of
$R$-modules $0\rightarrow P \rightarrow P_{n-1}\rightarrow \cdots
\rightarrow P_{0} \rightarrow  F \rightarrow 0$, where each $P_i$ is
projective. Decomposing this sequence into short exact sequences and
applying successively \cite[Exercise 8.2, page 223]{Rot} (since $F$
is flat), we obtain the following exact sequences of $T$-modules and
$R/M$-modules, respectively:
$$
\begin{array}{c}
  0\rightarrow P \otimes_R T\rightarrow P_{n-1}\otimes_R T\rightarrow \cdots
\rightarrow P_{0} \otimes_R T \rightarrow  F\otimes_R T \rightarrow 0 \\
  0\rightarrow P/MP \rightarrow P_{n-1}/MP_{n-1}\rightarrow \cdots
\rightarrow P_{0}/M P_{0} \rightarrow  F/MF \rightarrow 0 \\
\end{array}
$$
Since $\cot.D(T)\leq n$ and $\cot.D(D)\leq n$,  $ P \otimes_R S$ is
a projective $S$-module and  $ P/MP$ is a projective $R/M$-module.
Thus, $P$ is a projective $R$-module (by \cite[Theorem
5.1.1(1)]{Glaz}), and then $\pd_R(F)\leq n$. This means that $R$ is
$n$-perfect and therefore $\cot.D(R)\leq n$, which implies the
desired inequality.\cqfd

\end{section}\bigskip \bigskip


\begin{section}{Applications and Examples}

Now, we are ready to give examples showing that each of the
inequalities of Theorems \ref{thm-C-gldim-wdim} and
\ref{thm-C-FPD-FFD} can be strict. We begin with two examples
concerning Theorem  \ref{thm-C-gldim-wdim}. To this aim, we need the
following result.

\begin{prop}[\cite{Dobs}, Proposition 2.1]\label{pro-gldim-pullback}
Let  $V$ be a valuation domain of the form $K+M$, where $K$ is a
field and $M(\not= 0)$ is a maximal ideal of $V$. Let $D$ be a
proper subring of $K$ such that $\qf(D)=K$. For the pullback ring
$R=D+M$, we have:
\begin{enumerate}
    \item If $n=\gldim(V)$ and  $m=\gldim(D)$, then:
   $$ \gldim(R)= \left\{%
\begin{array}{ll}
    $n $  & \hbox{if  $\ n>m$ ;} \\
    $m$   & \hbox{if $\ m \geq n\ $  and   $\ \pd_D(K)<m$ ;} \\
   $m+1$  & \hbox{if $\ m\geq n\ $ and   $\ \pd_D(K)=m$ .} \\
\end{array}%
\right. $$
    \item $\wdim(R)=\wdim(D)$.
\end{enumerate}
\end{prop}

The following  example  shows a ring $R$ of finite global dimension
with the strict inequality $\cot.D(R)< \gldim(R)$  holds true.

\begin{exmp}\label{exmp1}
Let $\Z$ denote the ring of integers, and let $\Q$ denote the field
of rational numbers. Consider the discrete  valuation domain
$V=\Q[X]_{(X)}=\Q+XV$ (the localization of the polynomial ring
$\Q[X]$ at $(X)$). Then, for the pullback ring $R=\Z+XV$, we have:
$\gldim(R)=2$, $\wdim(R)=1$, and $\cot.D(R)=1$.
\end{exmp}
\proof By Proposition \ref{pro-gldim-pullback},
$\wdim(R)=\wdim(\Z)=1$ and $\gldim(R)=2$ (since $\pd_{\Z}(\Q)=1$).\\
From Theorem \ref{thm-pullback-C-gldim}, $\cot.D(R)=
\sup\{\cot.D(V),\cot.D(\Z)\}=1$ (since
$\gldim(V)=\gldim(\Z)=1$).\cqfd\bigskip

Note that, in Example \ref{exmp1} above,
$\gldim(R)=\cot.D(R)+\wdim(R)$. This means that the upper bound on
the global dimension given in Theorem \ref{thm-C-gldim-wdim} is the
best possible upper bound.\bigskip

In the next example we provide a ring $R$ of finite global dimension
with the strict inequalities $\cot.D(R)< \gldim(R)< \cot.D(R) +
\wdim(R)$  hold true.

\begin{exmp}\label{exmp2} Let $\C$ denote
the field of complex numbers and $\C(X,Y)$ denote the quotient field
of the polynomial ring $\C[X,Y]$.  Consider the discrete valuation
domain  $V=\C(X,Y)[Z]_{(Z)}=\C(X,Y)+ZV$ (the localization of the
polynomial ring $\C(X,Y)[Z]$ at $(Z)$). Then, for the pullback ring
$R=\C[X,Y]+ZV$, we have:   $\gldim(R)=3$, $\wdim(R)=2$, and
$\cot.D(R)=2$.
\end{exmp}
\proof By Proposition \ref{pro-gldim-pullback},
$\wdim(R)=\wdim(\C[X,Y])=2$.\\
We have $2=\gldim(\C[X,Y])>\gldim(V)=1$  and
$\pd_{\C[X,Y]}(\C(X,Y))=2)$ (by
\cite[Theorem 2]{Kap}). Then, from Proposition \ref{pro-gldim-pullback},  $\gldim(R)=2+1=3$.\\
Since, $\cot.D(\C[X,Y])=2$ (by Theorem \ref{thm-Syzygy-C-gldim}), we
obtain, from Theorem \ref{thm-pullback-C-gldim}, that $\cot.D(R)=
\sup\{\cot.D(\C[X,Y]),\cot.D(V)\}=2$,  as desired.\cqfd

\begin{rem}\label{rem-dom}\textnormal{
Example \ref{exmp2} above  answers the question evoked in Section 3,
above Proposition \ref{pro-C-gldim-domain}.\\ Indeed, consider a
valuation domain of the form $V=K+M$ and $R=D+M$. Let $Q$ denote the
quotient field of $R$ and then of $V$. We have $\pd_R(Q)=\pd_V(Q)$
(this is obtained using the same argument as in the proof of the
second part of Theorem \ref{thm-pullback-C-gldim} and since
$Q\otimes_R V = Q$ and $Q\otimes_R R/M=0$).\\
Applying this fact to Example \ref{exmp2} above, we have
$\pd_R(Q)=\pd_V(Q)=1$, and so $1=\pd_R(Q)<\cot.D(R)=2$.}
\end{rem}

Now we set examples concerning Theorem  \ref{thm-C-FPD-FFD}. To this
aim, we use Theorem \ref{thm-product-C-gldim} and the fact that
$\gldim(S\times R)=\sup\{\gldim(S),\gldim(R)\}$ for two rings $R$
and $S$ \cite[Chapitre VI, Exercise 8, page 123]{CE}. This result is
a consequence of  \cite[Lemma 2.5(2)]{Mah2001}, which implies also
that $\FPD(S\times R)=\sup\{\FPD(S),\FPD(R)\}$. Similarly,
\cite[Lemma 3.7]{BM2} is used to prove that $\wdim(S\times
R)=\sup\{\wdim(S),\wdim(R)\}$ and $\FFD(S\times
R)=\sup\{\FFD(S),\FFD(R)\}$.\bigskip

Next example shows a ring $R$ which has finite finitistic projective
dimension, and infinite weak dimension and satisfies $\cot.D(R)<
\FPD(R)=\cot.D(R)+\FFD(R)$.

\begin{exmp}\label{exmp3}
Consider a perfect ring $S$ of infinite weak dimension. Then, for
the ring $R$ of Example \ref{exmp1}, we have: $\cot.D(S\times R)=1$,
$\FPD(S\times R)=2$,  $\FFD(S\times R)=1$, and $\wdim(S\times
R)=\infty$.
\end{exmp}
\proof Apply  the facts mentioned above.\cqfd\bigskip

Using again the facts mentioned above, we get an example of a ring
$R$ with the strict inequalities $\cot.D(R)< \FPD(R)<  \cot.D(R) +
\FFD(R)$ hold true.

\begin{exmp}\label{exmp5} Consider a perfect ring $S$ of infinite weak
dimension. Then, for the ring $R$ of Example \ref{exmp2}, we have:
$\cot.D(S\times R)=2$, $\FPD(S\times R)=3$,  $\FFD(S\times R)=2$,
and $\wdim(S\times R)=\infty$.
\end{exmp}
\bigskip
\noindent {\bf ACKNOWLEDGEMENTS.} The authors would like to express their sincere thanks for 
the referee for his/her helpful suggestions. \\

\end{section}


\end{document}